\documentclass[12pt]{article}
\usepackage{amsmath,amsfonts,amssymb,amsthm}
\newtheorem{theo}{Theorem}
\newtheorem*{conj}{Conjecture}
\newtheorem*{prop}{Lemma}

\title{A new family of positive integers}
\author{Michel Lassalle\\
\small Centre National de la Recherche Scientifique\\[-0.8ex]
\small Institut Gaspard Monge, Universit\'e de Marne-la-Vall\'ee\\[-0.8ex]
\small 77454 Marne-la-Vall\'ee Cedex, France\\[-0.8ex]
\small \texttt{lassalle @ univ-mlv.fr}\\[-0.8ex]
\small \texttt{http://www-igm.univ-mlv.fr/{\textasciitilde}lassalle/index.html}}
\date{\small 2000 Mathematics Subject Classification : 05A10, 33C20}

\begin{document}
\maketitle
\begin{abstract}
Let $n,p,k$ be three positive integers. We prove that the numbers
$\binom{n}{k} \ {}_{3}^{}{F}_{2} (1-k,-p,p-n\ ; \ 1,1-n \ ; \ 1)$ 
are positive integers which generalize the classical binomial coefficients.
We give two generating functions for these integers, and a 
straightforward application.
\end{abstract}

\section{Definition}
We use the standard notation 
for hypergeometric series
\[{}_pF_q\!\left[\begin{matrix}a_1,a_2,\dots,a_p\\
b_1,b_2,\dots,b_q\end{matrix};z\right]=
\sum _{k=0} ^{\infty}\frac {(a_1)_k\dots(a_p)_k}
{(b_1)_k\dots(b_q)_k}\frac{z^k}{k!},\]
where for an indeterminate $a$ and some positive integer $k$,  
the raising factorial is defined by ${(a)}_k = a(a+1) \ldots (a+k-1)$. 

There are not many families of positive integers which may be 
defined in terms of hypergeometric functions. Among them stand of 
course the binomial coefficients
\[\binom{n}{p} = {}_2F_1\!\left[\begin{matrix}-p,p-n\\
1\end{matrix};1\right].\]
Actually this expression is obtained by specializing $x=p-n, y=1$ in the 
celebrated Chu - Vandermonde formula
\[{\frac{{(y-x)}_p}{{(y)}_p}} = {}_2F_1\!\left[\begin{matrix}-p,x\\
y\end{matrix};1\right].\]

Of course using this  relation as a definition of binomial coefficients would 
be rather tautological. However, quite surprisingly, it is possible to define a 
new family of positive integers by slightly modifying the 
Chu - Vandermonde formula.

Indeed for any positive integers $n,p,k$, let us define
\[{\binom{n}{p}}_{k} = \binom{n}{k} \ 
{}_3F_2\!\left[\begin{matrix}1-k,-p,p-n\\
1-n,1\end{matrix};1\right].\]
We have obviously
\[{\binom{n}{p}}_{k} = 0 \quad \textrm{for} \quad k>n \quad , \quad
{\binom{n}{p}}_{k} ={\binom{n}{n-p}}_{k} \quad , \quad 
{\binom{n}{p}}_{n} =\binom{n}{p},\]
the last equation following directly from the Chu - Vandermonde formula.

This definition can be rewritten
\begin{equation}
{\binom{n}{p}}_{k} = \binom{n}{k} \ \sum_{r \ge 0} 
\binom{p}{r}\binom{n-p}{r} \frac{\binom{k-1}{r}}{\binom{n-1}{r}}
= \frac{n}{k} \ \sum_{r \ge 0} 
\binom{p}{r}\binom{n-p}{r} \binom{n-r-1}{k-r-1}.
\end{equation}
Thus $\frac{k}{n}\ {\binom{n}{p}}_{k}$ is a positive integer. 

One has easily
\begin{equation*}
{\binom{n}{0}}_{k} = \binom{n}{k} \quad, \quad
{\binom{n}{1}}_{k} = k \binom{n}{k}\quad, \quad
{\binom{n}{2}}_{k} = k \binom{n}{k} 
+\frac{n(n-3)}{2} \binom{n-2}{k-2},
\end{equation*} 
\[{\binom{n}{p}}_{0} = 0 \quad \textrm{for} \quad p \neq 0,n \quad , 
\quad {\binom{n}{p}}_{1} = n,
\quad {\binom{n}{p}}_{2} = \frac{n}{2} \left(n-1+p(n-p)\right),\]
and also
\begin{equation*}
\begin{split}
{\binom{n}{p}}_{n-1} &= n\left[\binom{n-1}{p-1} 
+\binom{n-2}{p}\right],\\
{\binom{n}{p}}_{n-2}
&= \binom{n}{2}\left[\binom{n-2}{p}+\binom{n-2}{p-2}\right]
+\frac{n(n-3)}{2}\binom{n-4}{p-2}.
\end{split}
\end{equation*}
These relations suggest that ${\binom{n}{p}}_{k}$ is a positive integer.

\section {Integrality}

Using the Chu - Vandermonde formula
\[\binom{n-r-1}{n-k} = \sum_{i = 0}^{n-k} {(-1)}^i \binom{r}{i} 
\binom{n-i-1}{n-k-i},\]
we get
\[{\binom{n}{p}}_{k} =\frac{n}{k} \ \sum_{i = 0}^{n-k} 
{(-1)}^i \binom{n-i-1}{k-1} \sum_{r \ge 0}
\binom{p}{r} \binom{n-p}{r} \binom {r}{i}.\]
But using again the Chu - Vandermonde formula one has 
\begin{equation*}
\sum_{r \ge 0} \binom {r}{i}
\binom{p}{r} \binom{n-p}{r} = \binom{p}{i} 
\sum_{r \ge 0} \binom {p-i}{r-i} \binom {n-p}{r}
= \binom{p}{i} \binom {n-i}{p}.
\end{equation*}
Hence we obtain
\begin{equation}
\begin{split}
{\binom{n}{p}}_{k} &=\frac{n}{k} \ \sum_{i = 0}^{n-k} 
{(-1)}^i \binom{n-i-1}{k-1} \binom{p}{i} \binom {n-i}{p}\\
&=\frac{n}{k} \ \sum_{i = 0}^{n-k} 
{(-1)}^i \binom{n-i-1}{k-1} \binom{n-i}{i} \binom {n-2i}{p-i}\\
&= \sum_{i = 0}^{n-k} {(-1)}^i \binom{n-i}{k} \frac{n}{n-i} 
\binom{n-i}{i} \binom {n-2i}{p-i}.
\end{split}
\end{equation}

Finally we get
\[{\binom{n}{p}}_{k}=
\sum_{i = 0}^{n-k} {(-1)}^i \binom{n-i}{k} \binom {n-2i}{p-i}
\left[\binom {n-i}{i} + \binom {n-i-1}{i-1}\right].\]
We have thus proved
\begin{theo}
The positive number ${\binom{n}{p}}_{k}$ is an integer.
\end{theo}

Thanks are due to Jiang Zeng for shortening the proof of this result.
Note that the previous relations imply immediatly

\begin{equation*}
\begin{split}
{\binom{n}{p}}_{k}+{\binom{n}{p-1}}_{k} &=
\sum_{i = 0}^{n-k} {(-1)}^i  \binom{n-i}{k} \binom {n-2i+1}{p-i}
\left[\binom {n-i}{i} + \binom {n-i-1}{i-1}\right]\\
&=\frac{k+1}{n+1} \ {\binom{n+1}{p}}_{k+1} -
\sum_{i = 2}^{n-k} {(-1)}^i \binom{n-i}{k} \binom {n-2i+1}{p-i} 
\binom{n-i-1}{i-2}.
\end{split}
\end {equation*}

An intriguing problem is to get a combinatorial interpretation for 
${\binom{n}{p}}_{k}$.

\section{Generating functions}

The following generating function is due to Jiang Zeng.

\begin{theo}
We have
\begin{multline*}
\sum_{k,p \geq 0} {\binom{n}{p}}_{k} x^p y^k = {2}^{-n}
\Big[ \left( (1+x)(1+y)+\sqrt{(1+x)^2(1+y)^2-4x(1+y)} \right)^n \\
+\left((1+x)(1+y)-\sqrt{(1+x)^2(1+y)^2-4x(1+y)}\right)^n \Big].
\end{multline*}
\end{theo}

\begin{proof}
From equation (2) we get 
\begin{equation*}
\begin{split}
\sum_{k, p \geq 0}{\binom{n}{p}}_{k}\ x^p y^k&= \sum_{k, p \geq 0}
x^p y^k \sum_{i \geq 0}(-1)^i \binom{n-i}{k} \frac{n}{n-i} 
\binom{n-i}{i} \binom {n-2i}{p-i}\\
&= \sum_{i \geq 0}
(-1)^i (1+y)^{n-i} \frac{n}{n-i} \binom{n-i}{i} x^i(1+x)^{n-2i}\\
&= (1+x)^n(1+y)^n\sum_{i\geq 0}\frac{n}{n-i}\binom{n-i}{i}z^i,
\end{split}
\end{equation*}
with $z=-x/((1+x)^2(1+y))$. But we have the following
identity
\begin{equation*}
\sum_{n > i}\frac{n}{n-i}\binom{n-i}{i}z^i=\left(\frac{1+\sqrt{1+4z}}{2}\right)^n
+\left(\frac{1-\sqrt{1+4z}}{2}\right)^n.
\end{equation*}
\end{proof}

We can give another generating function. The following recurrence 
relation is needed.
\begin{prop}
We have
\[(n-p+1){\binom{n}{p-1}}_k-p{\binom{n}{p}}_k
=\frac{n}{n-1}(n-2p+1){\binom{n-1}{p-1}}_k.\]
\end{prop}
\begin{proof}
This can be easily deduced from equation (1).
Indeed up to $n/k$ the left-hand side can be written 
\begin {equation*}
\begin{split}
&\sum_{r \ge 0} \binom{n-r-1}{k-r-1}
\left[ (n-p+1) \binom{p-1}{r}\binom{n-p+1}{r} -p 
\binom{p}{r}\binom{n-p}{r} \right]\\  &=
(n-2p+1)\sum_{r \ge 0} \binom{n-r-1}{k-r-1}
\left[\binom{p-1}{r} \binom{n-p}{r} - \binom{p-1}{r-1} 
\binom{n-p}{r-1}\right] \\&=
(n-2p+1)\sum_{r \ge 0} \binom{p-1}{r} \binom{n-p}{r}
\left[\binom{n-r-1}{k-r-1}
- \binom{n-r-2}{k-r-2} \right]\\ &=
(n-2p+1)\sum_{r \ge 0} 
\binom{p-1}{r} \binom{n-p}{r}\binom{n-r-2}{k-r-1}.
\end{split}
\end {equation*}
\end{proof}

\begin{theo}
We have
\[\sum_{k \ge 1} {\binom{n}{p}}_{k} y^k =
 n y (y+1)^p \,{}_2F_1\!\left[\begin{matrix}p+1,p-n+1\\
2\end{matrix};-y\right].\]
\end{theo}

\begin{proof} By recurrence over the integers $n$ and $p$. 
The property is true for $p=0$, since we have
\[\sum_{k \ge 1} {\binom{n}{0}}_{k} y^k =
\sum_{k \ge 1} \binom{n}{k} y^k = (1+y)^n-1 =
 n y  \,{}_2F_1\!\left[\begin{matrix}1,1-n\\
2\end{matrix};-y\right].\]
From the previous recurrence relation, we deduce that it is enough to 
prove
\begin{multline*}
(n-p+1) \ {}_2F_1\!\left[\begin{matrix}p,p-n\\
2\end{matrix};-y\right]
-p (1+y) \ {}_2F_1\!\left[\begin{matrix}p+1,p-n+1\\
2\end{matrix};-y\right]=\\
(n-2p+1) \ {}_2F_1\!\left[\begin{matrix}p,p-n+1\\
2\end{matrix};-y\right].
\end{multline*}
But this is a classical contiguity relation for ${}_2F_1$
(see for instance~\cite{Ra}, Exercice 21.8, page 71), namely
\[(c-b-1) \ {}_2F_1\!\left[\begin{matrix}a,b\\
c\end{matrix};y\right]
-a (1-y) \ {}_2F_1\!\left[\begin{matrix}a+1,b+1\\
c\end{matrix};y\right]
=(c-a-b-1) \ {}_2F_1\!\left[\begin{matrix}a,b+1\\
c\end{matrix};y\right].\]
\end{proof}

By identification of the coefficients of $y$ we obtain a property 
which seems difficult to be proved directly, namely
\begin{equation*}
\begin{split}
{\binom{n}{p}}_{k}&= \frac{n}{p}
 \sum_{i=0}^{k-1} \binom{n-p+i}{n-p} \binom{p}{i+1}
\binom{n-p}{k-i-1}\\ &= \frac{n}{p}
 \sum_{i=0}^{n-k} \binom{k-1+i}{n-p} \binom{p}{n-k-i} \binom{n-p}{i}.
\end{split}
\end{equation*}
Hence $\frac{p}{n} {\binom{n}{p}}_{k}$ is a positive integer. 

This gives a second proof that ${\binom{n}{p}}_{k}$ is a positive integer. 
Indeed we get
\begin{equation*}
\begin{split}
{\binom{n}{p}}_{k} &= \frac{p}{n} {\binom{n}{p}}_{k} +\frac{n-p}{n} 
{\binom{n}{n-p}}_{k}\\&=
\sum_{i=0}^{n-k} \left[ \binom{k-1+i}{n-p} \binom{p}{n-k-i} \binom{n-p}{i}
+ \binom{k-1+i}{p} \binom{n-p}{n-k-i} \binom{p}{i} \right].
\end{split}
\end{equation*}
Incidentally we have also proved
\begin{multline*}
{2}^{-n}
\Big[ \left( (1+x)(1+y)+\sqrt{(1+x)^2(1+y)^2-4x(1+y)} \right)^n \\
+\left((1+x)(1+y)-\sqrt{(1+x)^2(1+y)^2-4x(1+y)}\right)^n \Big] =\\ 
1+ x^n + n y \sum_{p=0}^n x^p (y+1)^p 
\,{}_2F_1\!\left[\begin{matrix}p+1,p-n+1\\2\end{matrix};-y\right].
\end{multline*}

\section {Theory of partitions}

Let us now indicate in which situation the integers ${\binom{n}{p}}_{k}$ 
are naturally encountered. 

A partition $\lambda= ( {\lambda }_{1},...,{\lambda }_{n})$ 
is a finite weakly decreasing
sequence of positive integers, called parts. The number 
$n=l(\lambda)$ of parts is called the length of 
$\lambda$, and
$|\lambda| = \sum_{i = 1}^{n} \lambda_{i}$
the weight of $\lambda$. For any integer $i\geq1$, 
${m}_{i} (\lambda)  = \textrm{card} \{j: {\lambda }_{j}  = i\}$
is the multiplicity of $i$ in $\lambda$.  We set 
\[{z}_{\lambda }  = \prod\limits_{i \ge  1}^{} 
{i}^{{m}_{i}(\lambda)} {m}_{i}(\lambda) !  .\]

Let $X$ be an indeterminate and $n$ a positive integer. We write 
${[X]}_{n }  =  X (X - 1) \cdots (X -n +1)$	
for the lowering factorial and
$\binom{X}{n}  =  {[X]}_{n}/n!$. 
The following result has been proved in~\cite{La}
(Theorem 4, p. 275) and in~\cite{Z} : for any positive integers $n,s$ we have
\begin{equation*}
\sum_{\left|{\mu }\right| = n} 
{\frac{{X}^{l(\mu ) - 1}}
{{z}_{\mu }}}  \left({ \sum_{i = 1}^{l(\mu )} 
{({\mu }_{i})}_{s} }\right)
= (s - 1) !  \sum_{k = 1}^{\min (n,s)}  
\binom{s}{k} \binom{X +n - 1}{n - k} .
\end{equation*}
This property generalizes as follows.
\begin{theo} 
Let $X$ be an indeterminate and $n,r,s$ three positive integers. We have
\begin{equation*}
\sum_{\left|{\mu }\right| = n} 
{\frac{{X}^{l(\mu ) - 1}}
{{z}_{\mu }}}  \left( \sum_{i = 1}^{l(\mu )} 
{({\mu }_{i})}_{r} {({\mu }_{i})}_{s} \right)
= \frac{r! s!}{r+s}  \sum_{k = 1}^{\min (n,r+s)}  
{\binom{r+s}{s}}_k \binom{X +n - 1}{n - k} .
\end{equation*}
\end{theo}

\begin{proof}
By recurrence over $r$. For $r=0$  the property has been proved since
${\binom{s}{s}}_k=\binom{s}{k}$. Now one has 
$(i)_{r+1}(i)_{s}=(i)_{r}(i)_{s+1}+(r-s)(i)_r(i)_s$. Thus
it is enough to prove
\[\frac{(r+1)! s!}{r+s+1}{\binom{r+s+1}{s}}_k=
\frac{r! (s+1)!}{r+s+1}{\binom{r+s+1}{s+1}}_k+
(r-s)\frac{r! s!}{r+s}{\binom{r+s}{s}}_k.\]
But this is the statement of the Lemma.
\end{proof}

This result suggests the following conjecture.

\begin{conj}
Let $X$ be an indeterminate, $n$ a positive integer and 
$r=(r_1,\ldots,r_m)$ a positive multi-integer with weight
$|r| = \sum_{i = 1}^m r_i$. We have
\begin{equation*}
\sum_{|\mu| = n} 
\frac{X^{l(\mu) - 1}}
{z_{\mu }}  \left( \sum_{i = 1}^{l(\mu )} 
\prod_{k = 1}^m {(\mu_i)}_{r_k} \right)
= \frac{\prod_{j} r_{j}!}{|r|}  \sum_{k = 1}^{\min (n,|r|)}  
c_k^{(r)} \binom{X +n - 1}{n - k} ,
\end{equation*}
where the coefficients $c_k^{(r)}$ are positive integers, to be 
computed.
\end{conj}

The integers ${\binom{n}{p}}_{k}$ appear also when studying (shifted) Jack 
polynomials in the spirit of~\cite{La1,La2} 
(this application will be given in a forthcoming paper).

\end{document}